\documentclass{commatDV}

\theorembodyfont\it
 \newtheorem{tthr}[theorem]{Theorem}
 \newtheorem{corr}[theorem]{Corollary}
 \newtheorem{lmr}[theorem]{Lemma}
 \newtheorem{prr}[theorem]{Proposition}
 \newtheorem{hyp}[theorem]{Conjecture}

\theorembodyfont\rm
 \newtheorem{defnr}{Definition}
 \newtheorem{remnr}{Remark}

\newcommand{\btr}{
\begin{tthr}
}
\newcommand{\etr}{
\end{tthr}
}
\newcommand{\bpr}{
\begin{prr}
}
\newcommand{\epr}{
\end{prr}
}
\newcommand{\bcr}{
\begin{corr}
}
\newcommand{\ecr}{
\end{corr}
}
\newcommand{\blr}{
\begin{lmr}
}
\newcommand{\elr}{
\end{lmr}
}
\newcommand{\bh}{
\begin{hyp}
}
\newcommand{\eh}{
\end{hyp}
}
\newcommand{\bdr}{
\begin{defnr}
}
\newcommand{\edr}{
\end{defnr}
}
\newcommand{\brm}{
\begin{remnr}
}
\newcommand{\erm}{
\end{remnr}
}

\newcommand{\bp}{
\begin{proof}
}
\newcommand{\ep}{
\end{proof}
}
\newcommand{\beq}{
\begin{equation}
}
\newcommand{\eeq}{
\end{equation}
}
\newcommand{\bcs}{
\begin{cases}
}
\newcommand{\ecs}{
\end{cases}
}
\newcommand{\en}{
\begin{enumerate}
}
\newcommand{\ene}{
\end{enumerate}
}

\newcommand{\D}{\Delta}
\newcommand{\lab}{\label}
\def\>{\geqslant}
\def\<{\leqslant}

\newcommand{\quat}[1]{\textquotedblleft #1\textquotedblright}
\newcommand{\mpp}[1]{\;\mathrm{mod}\,#1}

\title{On a theorem of J.~Shallit concerning Fibonacci partitions}
\author{F.~V.~Weinstein}
\affiliation{N/A
    \email{felix.weinstein46@gmail.com}%
    }

\abstract{%
In this note, I prove a~claim on determinants of some special tridiagonal matrices. Together with my result about Fibonacci partitions (\url{https://arxiv.org/pdf/math/0307150.pdf}), this claim allows one to prove one (slightly strengthened) Shallit's result (\url{https://arxiv.org/pdf/2007.14930.pdf}) about such partitions.
}



\VOLUME{30}
\NUMBER{3}
\YEAR{2022}
\firstpage{203}
\DOI{https://doi.org/10.46298/cm.10769}

\begin{paper}

\section{Introduction}\lab{Intr}

Let $f_1 = 1,f_2 = 2$ and $f_i = f_{i-1}+f_{i-2}$ for $i>2$ be the sequence of Fibonacci numbers. Observe that the \quat{conventional} definition of Fibonacci numbers is different, see \url{http://en.wikipedia.org/wiki/Fibonacci_number}.

A \textit{Fibonacci partition} of a~positive integer
$n$ is a~representation of $n$ as an unordered sum of distinct
Fibonacci numbers, which are referred to as the \textit{parts} of the Fibonacci partition.

Let $\Phi_h(n)$ be the quantity (the cardinality of the set) of Fibonacci partitions of $n$ with $h$ parts. J.~Shallit has established the following interesting property of the function $\Phi_h(n)$:
for integers $n>0$, ${d\>2}$ and $i$, let $r_{d,i}(n)$ be the quantity of all Fibonacci partitions of $n$
with number of parts $\equiv i\mpp d$. Then, (see
\tracingmacros = 1~\cite[Th.~2]{SH})
\tracingmacros = 0
\[
\left|r_{3,i}(n)-r_{3,i+1}(n)\right|\<1.
\]
To prove this inequality, J.~Shallit used a~technique of automata theory.

Set
\[
\Phi(n;t): = \sum_{h>0}\Phi_h(n)t^h.
\]
In~\cite{FV}, I obtained a~formula
which expresses $\Phi(n;t)$ as determinant of a~tridiagonal matrix depending on $n$.
In \S\ref{Sp} of this note, I establish Theorem~\ref{t1}
on a~property of such determinants.

In \S\ref{Fb}, I explain
(Theorem~\ref{t2}) how the mentioned (see~\cite{FV}) formula for $\Phi(n;t)$ together with Theorem~\ref{t1}
imply not only Shallit's result, but also the formula
\[
\big(r_{3,0}(n)-r_{3,1}(n)\big)\cdot\big(r_{3,0}(n)-r_{3,2}(n)\big)\cdot\big(r_{3,1}(n)-r_{3,2}(n)\big) = 0.
\]

\section{3--special polynomials}\lab{Sp}

Let $d\>2$ be an integer number. For any $g(t) = \sum_{h\>0}a_ht^h \in\mathbb{Z}[t]$, define
\[
\|g(t)\|: = \sum_{h\>0}\;a_h,\ \ R_i(g(t)): =
\sum_{h\equiv i\mpp d}a_h,\ \text{where\; $i\in\{0,1,\dots,d-1\}$}.
\]

Let $K_d[T]: = \mathbb{Z}[T]/(T^d -1)$.
Define a~map $R^{(d)}:\mathbb{Z}[t]\rightarrow K_d[T]$ by the formula
\[
R^{(d)}(g(t)): = R_0(g(t))+R_1(g(t))T+\dots+R_{d-1}(g(t))T^{d-1}.
\]

The following Lemma is subject to easy direct verification.
\blr\lab{lm1}
The map $R^{(d)}:\mathbb{Z}[t]\rightarrow K_d[T]$ is a~homomorphism of $\mathbb{Z}$--algebras.
\elr

In this Section, I consider only the case $d = 3$. For brevity, set $K: = K_3[T]$ and $R: = R^{(3)}$.

For any $g(t)\in\mathbb{Z}[t]$, we obviously have
\beq\lab{eq0}
R\big((1+t+t^2)\cdot g(t)\big) = \|g(t)\|\cdot\varphi(T),\text{\quad where $\varphi(T): = 1+T+T^2$.}
\eeq

\bdr
We say that $a+b\,T+c\,T^2 \in K$ is a~\emph{special element} if either $a = b = c$, or
 $|a-b|+|a-c|+|b-c| = 2$.
\edr

Formula~\eqref{eq0} easily implies

\blr\lab{lmm}
An element $A[T]\in K$ is special if and only if
\[
A[T]\cdot(T-1)\in M[T]: = \left\{0,\pm\;(T-1),\pm\;T(T-1),\pm\;T^2 (T-1)\right\}.
\]
\elr

\bcr\lab{lm3}
Any product of special elements is a~special element.
\ecr

\bdr
We say that $g(t)\in\mathbb{Z}[t]$ is a~3--\emph{special polynomial} if $R(g(t))$ is a~special element.
\edr

In what follows, $A = (a_1,a_2,\dots,a_m)$ is either a~vector with integer non-negative coordinates if $m>0$,
or the empty set if $m = 0$.
Let us define a~polynomial
\[
\D(A;t): = \D(a_1,\dots,a_m;t)\in\mathbb{Z}[t]
\]
by the formulas
\[
\D(\emptyset;t): = 1,\qquad\D(0;t): = 0,\qquad\D(a;t): = t+t^2 +\dots+t^a \quad\text{for $a>0$},
\]
\beq\lab{eq4}
\D(a_1,\dots,a_m;t): =
\D(a_1,\dots,a_{m-1};t)\cdot\D(a_m;t)
-\D(a_1,\dots,a_{m-2};t)\cdot t^{a_m+1} \quad \text{if $m\>2$.}
\eeq
 Obviously, for $m>0$,
\[
\footnotesize
\D\left(a_1,a_2,\dots,a_m;t\right) =
\begin{vmatrix}
\;\D(a_1;t) & t^{a_2+1} & 0 & 0 & \hdots & 0 \\
\;1 & \D(a_2;t) & t^{a_3+1} & 0 & \hdots & 0 \\
\;\vdots & \vdots & & \ddots & &\vdots \\
\;0 & 0 & \hdots & 1 &\D(a_{m-1};t) & t^{a_m+1} \\
\;0 & 0 & \hdots & 0 & 1 &\D(a_m;t)
\end{vmatrix}
.
\]

The main result of this note is
\btr\lab{t1}
For any $A = (a_1,a_2,\dots,a_m)$, the polynomial $\D(A;t)$ is a~\textrm{3}--special one.
\etr
The proof uses the following auxiliary claim.

\blr\lab{lm4}
Let $\varepsilon(A): = (\varepsilon(a_1),\dots,\varepsilon(a_m))$,
where $\varepsilon(a): = a-3\left\lfloor\frac{a}{3}\right\rfloor$.
Then,
\[
R\big(\D(A;t)\big) = R\left(\D\big(\varepsilon(A);t\big)\right)
+k\cdot\varphi(T),
\qquad\text{where\; $k = k(A)\in\mathbb{Z}$}.
\]
\elr

\bp
Let us prove by induction on $m$. For $m = 1$ and $a\>1$, we have
\[
\D(a;t) = t\left(1+t^3 +\dots+t^{3\left\lfloor\frac{a}{3}\right\rfloor}\right)(1+t+t^2)+
t^{3\left\lfloor\frac{a}{3}\right\rfloor}\cdot\D\big(\varepsilon(a);t\big).
\]
Applying $R$ to both sides of this equality we obtain
\beq\lab{eq5}
R\big(\D(a;t)\big) = R\big(\D\big(\varepsilon(a);t\big)\big)+k\cdot\varphi(T),
\qquad\text{where\; $k = 1+\left\lfloor\frac{a}{3}\right\rfloor$}.
\eeq
\smallskip
For $m\>2$, let us apply $R$ to expression~\eqref{eq4}.
The induction hypothesis, Lemma~\ref{lm1},
formulas~\eqref{eq0} and~\eqref{eq5}, the obvious formula $R(t^a) = T^{\varepsilon(a)}$,
and a~short computation yield the required result.
\ep

\begin{Proof}
{Proof of Theorem $\ref{t1}$}
In view of Lemma~\ref{lm4}, it suffices to assume that
$a_i\in\{0,1,2\}$ for any $i = 1,2,\dots,m$.
Keeping Lemma~\ref{lmm} in mind, define
\[
S(a_1,\dots,a_m): = R\big(\D(a_1,\dots,a_m;t)\big)\cdot(T-1)\in K.
\]
The expression~\eqref{eq4} and formula $\varphi(T)\cdot(T-1) = 0$
easily imply the recurrent formula
\beq\lab{eqs}
\footnotesize
S(a_1,\dots,a_m) =
\bcs
-S(a_1,\dots,a_{m-2})\cdot T&\text{if $a_m = 0$},\\
S(a_1,\dots,a_{m-1})\cdot T+S(a_1,\dots,a_{m-2})\cdot(T+1)&\text{if $a_m = 1$},\\
-S(a_1,\dots,a_{m-1})-S(a_1,\dots,a_{m-2})&\text{if $a_m = 2$}.
\ecs
\eeq
By Lemma~\ref{lmm} it remains to show that $S(a_1,\dots,a_m)\in M[T]$.

Let us prove this by induction on $m$.
For $m = 1,2$, the claim is directly checked.
In particular, $S(0) = 0$ and $S(a,0) = -T(T-1)$.

For $a_m = 0$, the last expressions and formula~\eqref{eqs} imply the theorem
by induction for any $m\>1$. Therefore, assume that $a_m = 1$ or $a_m = 2$.
From expressions~\eqref{eqs} it is not difficult to obtain the expressions
\[
S(a_1,\dots,a_{m-1},1) =
\bcs
S(a_1,\dots,a_{m-2},2)\cdot T^2 &\text{if $a_{m-1} = 0$},\\
-S(a_1,\dots,a_{m-2},2)\cdot(T+1)&\text{if $a_{m-1} = 1$},\\
S(a_1,\dots,a_{m-3},a_{m-2}+2)\cdot T&\text{if $a_{m-1} = 2$},
\ecs
\]

\[
S(a_1,\dots,a_{m-1},2) =
\bcs
-S(a_1,\dots,a_{m-3},a_{m-2}+2)\cdot T&\text{if $a_{m-1} = 0$},\\
S(a_1,\dots,a_{m-2},2)\cdot(T+1)&\text{if $a_{m-1} = 1$},\\
S(a_1,\dots,a_{m-3})&\text{if $a_{m-1} = 2$}.
\ecs
\]
Since
\[
(T-1)(T+1) = -T^2 (T-1),
\]
these expressions and the induction
hypothesis complete the proof.
\end{Proof}
\section{An application to Fibonacci partitions}\lab{Fb}

In \S 2 of the article~\cite{FV}, for any positive integer $n$,
a~certain sequence is uniquely defined
\beq\lab{eqa}
\alpha(n) = \left\{\alpha_1(n),\alpha_2(n),\dots,\alpha_{k(n)}(n)\right\}
\eeq
where $\alpha_k(n)$
is a~vector with positive integer coordinates for any $k = 1,2,\dots,k(n)$,
 and it is shown (\cite[Th.2.11]{FV}) that
\[
\Phi(n;t) = \D\big(\alpha_1(n);t\big)\cdot\D\big(\alpha_2(n);t\big)\cdot{\dots}\cdot\D\big(\alpha_{k(n)}(n);t\big).
\]
By Theorem~\ref{t1} the polynomial $\D\big(\alpha_k(n);t\big)$ is a~3--special one for any $k$.
Thus, Lemma~\ref{lm1} and Corollary~\ref{lm3} imply
\btr\lab{t2}
For any integer $n>0$, the polynomial $\Phi(n;t)$ is a~\textrm{3}--special one.
\etr

\brm
Using arguments similar to those in \S\ref{Sp} (where $d = 3$ is replaced with $d = 2$) and
the formula for $\Phi(n;t)$ one can easily show that
$|r_{2,0}(n)-r_{2,1}(n)|\<1$ for any positive integer $n$.
It is obvious that this inequality is equivalent to the analytic identity
\beq\lab{Id}
\prod_{i = 1}^{\infty}(1-x^{f_i}) =
1+\sum_{n = 1}^{\infty}\chi(n)x^n,\qquad\text{where}\quad|\chi(n)|\<1.
\eeq
For other proofs of this identity, see~\cite{AR},\cite{RB} and~\cite{FV}.

In addition to that, an interesting result of Y.~Zhao should be mentioned.
Namely, Proposition~2 of the article~\cite{Zh} implies the polynomial identity
\[
\prod_{a\<i\<b}(1-x^{f_i}) =
1+\sum_{n}\chi_{a,b}(n)x^n,\qquad\text{where}\quad\left|\chi_{a,b}(n)\right|\<1,
\]
which is valid for any positive integers $a\<b$.
\erm

\bh
For positive integers $a\<b$, let
\[
M(a,b): = \{f_a,f_{a+1},\dots,f_b\}.
\]
For integers $n>0$ and $i$,
let $r_{3,i}^{(a,b)}(n)$ be the quantity of Fibonacci partitions of $n$
with parts from the set $M(a,b)$ and
with number of parts $\equiv i\mpp 3$.

Then, $\left|r_{3,i}^{(a,b)}(n)-r_{3,j}^{(a,b)}(n)\right|\<1$
for any $i,j\in\{0,1,2\}$. Moreover,
\[
\left(r_{3,0}^{(a,b)}(n)-r_{3,1}^{(a,b)}(n)\right)
\cdot\left(r_{3,0}^{(a,b)}(n)-r_{3,2}^{(a,b)}(n)\right)
\cdot\left(r_{3,1}^{(a,b)}(n)-r_{3,2}^{(a,b)}(n)\right) = 0.
\]
\eh

\end{paper}